\documentclass[12pt]{article}

\usepackage{enumitem}

\usepackage[nolists]{endfloat}

\usepackage{mathrsfs}
\usepackage{mathptmx}       
\usepackage{helvet}         
\usepackage{courier}        
\usepackage{type1cm}        
\usepackage[numbers]{natbib}       
\usepackage{makeidx}         
\usepackage{graphicx}        
\usepackage{subfigure}     
\usepackage{rotating}            
\usepackage{multicol}        
\usepackage[bottom]{footmisc}
\usepackage{titlesec}
\usepackage{amsfonts}
\usepackage{amsmath,amssymb}
\usepackage{subfig}
\usepackage{url}
\usepackage{bm}
\usepackage{lineno}
\usepackage{floatrow,fr-subfig}
\usepackage[nolists]{endfloat}
\usepackage{rotating}

\newtheorem{theorem}{Theorem}
\newtheorem{corollary}{Corollary}

\newtheorem{proposition}{Proposition}
\newtheorem{assumption}{A}

\makeatletter
\renewcommand\section{\@startsection{section}{1}{\z@}%
                                      {-3.5ex \@plus -1ex \@minus -.2ex}%
                                      {2.3ex \@plus.2ex}%
                                      {\normalfont\bfseries}}
\makeatother
\makeatletter
\renewcommand\subsection{\@startsection{subsection}{1}{\z@}%
                                      {-3.5ex \@plus -1ex \@minus -.2ex}%
                                      {2.3ex \@plus.2ex}%
                                      {\noindent\normalfont\emph}}
\makeatother

\begin{document}

\begin{center} \begin{Large} \textbf{The Effects of Adaptation on Maximum Likelihood Inference for Non-Linear Models with Normal Errors}\end{Large} \\ [18 pt]
Nancy Flournoy$^\star$, Caterina May$^\dagger$ and Chiara Tommasi$^\ddagger$\\ [12 pt]\begin{small}\textit{$^\star$University of Missouri, Columbia, USA\\
		$^\dagger$Universit\`a degli Studi del Piemonte Orientale, Italy\\
$^\ddagger$Universit\`a degli Studi di Milano, Italy} \end{small} \\ [18 pt]
\end{center}
\vspace{2cm}

\author{Nancy Flournoy, Caterina May and Chiara Tommasi}

\noindent \textbf{Abstract:} This work studies the properties of the maximum likelihood estimator (MLE) of a non-linear model with Gaussian errors and multidimensional parameter. The observations are collected in a two-stage experimental design and are dependent since the second stage design is determined by the observations at the first stage; the MLE maximizes the total likelihood.\\
Differently from the most of the literature, the first stage sample size is small, and hence asymptotic approximation is used only in the second stage.
It is proved that the MLE is consistent and that its asymptotic distribution is a specific Gaussian mixture, via stable convergence.
Finally, a simulation study is provided in the case of a dose-response Emax model.
\noindent

\vspace{12pt}

\section{Introduction}
This paper deals with the problems related to the estimation of a non-linear multi-parameter model with Gaussian errors.
Optimal experimental design approach improves the efficiency of the estimate. As well known in the literature (see for instance \cite{Atk2007}), when an optimal experimental design is used to estimate the parameter of a non-linear model, the optimal design depends on the unknown parameter. A possibility to tackle this problem is to use a locally optimal design, which is based on a guessed value for the parameter. If this guessed value is poorly chosen, however, the locally optimal design may be poor too. 

One common approach to solve this problem is to adopt a two-stage procedure (see for instance \cite{Drag:Fedo:Wu:Adap:2008}, \cite{Dette:2012}, \cite{Pronz}). At the first stage an initial design is applied to collect the first-stage responses which are used to estimate the unknown parameter. This is the so called interim analysis. To collect the second stage responses, a locally optimal design is determined  using the estimated parameter from the interim analysis. Finally, the maximum likelihood method is applied to estimate the vector parameter, employing the whole sample of data. 

Note that the first and the second stage observations are dependent; the classical approach assumes that both stages have large sample dimensions, and hence the asymptotic theory can be applied, as in \cite{Dette:2012} and \cite{Pronz}. This approach eliminates the dependency between stages, which is mathematically useful, but not realistic in many applications.  In real life problems, in fact,  the sample size of the interim analysis can be small. Therefore, in this work we assume that only the second stage sample size goes to infinity  while the first stage sample size is fixed, and hence the standard asymptotic behavoiur of the maximum likelihood estimator (MLE) does not maintain. The present study extends \cite{Lane:Flou:2012} and \cite{Lane:Yao:Flou:2013} which considered 
	a unidimensional parameter and the design at each stage to be a single point; in \cite{Lane:Yao:Flou:2013} it is also shown, via simulations, that fixing the first stage sample size improves the limiting approximation;  this is an additional 
reason of the importance of the results here obtained.

Under these assumptions, we prove the consistency of the MLE. Furthermore, we prove that the asymptotic distribution of the MLE is a specific normal mixture; this is obtained via stable convergence (for an overview on stable convergence theory see \cite{Hausler}).
In this context of dependent data, the inverse of the Fisher information matrix is not the asymptotic covariance matrix of the MLE. However,  we provide an analytical relation between these two quantities, which justifies the idea of using a function of the information matrix as an optimality criterion. Finally, we compare the proposed two-stage adaptive design with a locally optimal design through a simulation study under the Emax model.
This study points out that there exist scenarios where the adaptive procedure is superior, although the behaviour is not symmetric with respect to the nominal values of the parameter. A tentative theoretical justification is given, based on the anlaytical expression of the first order bias term of the first stage MLE.

The paper is organized as follows. Section \ref{Sect:backg} recalls the basic concept and introduces the model and the notation. Section \ref{Sect:twostage} describes the two-stage adaptive experimental procedure and provides the structure of the likelihood in this particular case. Section \ref{sec:AsyProp} contains the main theoretical results. Section \ref{Sect:simulations} presents an example with simulations. In Section \ref{Sect:concl} a summary with a few comments conclude the paper.

\section{Background and Notation}\label{Sect:backg}
Assume $n$ independent observations follow  the model
\begin{align}\label{eq:Model0}
y_{j} = \eta(x_{j},\boldsymbol{\theta}) + \varepsilon_{j}, \quad \varepsilon_{j} \sim \mathscr{n}\left( 0,\sigma^2 \right),\quad j=1,\ldots,n,
\end{align}
 where $y_{j}$ is the response of the unit $j$ treated under an experimental condition $x_j\in {\cal X}$ and $\eta(x_{j},\boldsymbol{\theta})$ is some possibly non-linear continuous mean function of $p+1$ parameters, $\boldsymbol{\theta}=(\theta_0,\ldots,\theta_p)$, with $\boldsymbol{\theta}\in\boldsymbol{\Theta}$, where $\boldsymbol{\Theta}$ is a compact set in ${\mathbb R}^{p+1}$. In general, several units may be treated under the same experimental conditions. An  experimental design is a finite discrete probability distribution over ${\cal X}$:
\begin{equation}\label{xi}
\xi=\left\{\!\!\!\begin{array}{ccc}
x_1 &\cdots & x_M \\
\omega_1 & \cdots & \omega_M
\end{array}
\!\!\!  \right\},
\end{equation}
where $x_m$ denotes the {$m$th experimental point, or treatment, that may be} used in the {study} and $\omega_m$ {is} the proportion of {experimental} units to be taken at  {that} point{; $\omega_m\geq 0$ with  $\sum_{m=1}^{M}\omega_m=1$,} $ m=1,\ldots,M$ and $M$ is finite.

 It is well known that a good design can substantially improve the inferential results in a statistical analysis.
For instance, if the inferential goal is point estimation of $\boldsymbol{\theta}$, then an optimal design {may be} chosen {to} maximize some functional $\Phi(\cdot)$ of the information matrix
\begin{equation}\label{M}
M(\xi;\boldsymbol{\theta})=\int_{\cal X} \nabla \eta(x,\boldsymbol{\theta})
\nabla \eta(x,\boldsymbol{\theta})^T d\xi(x),
\end{equation}
as $M(\xi;\boldsymbol{\theta})^{-1}$ is proportional to the asymptotic covariance matrix of the maximum likelihood estimator (MLE) (\cite{Kief:Gene:1974}). 
In other terms, an optimal design for precise estimation of $\boldsymbol{\theta}$ is
\begin{equation}
\xi^*(\boldsymbol{\theta})=\arg\max_{\xi\in{\Xi}}\Phi[M(\xi;\boldsymbol{\theta})],
\label{eq:opt_dis}
\end{equation}
where $\Xi$ is the set of all the finite discrete probability distributions on ${\cal X}$ (i.e. the set of all designs). 
{Typically, taking derivatives  to approximate $\xi^*(\boldsymbol{\theta})$ will result in $\{\omega_m\}$
that are not multiples of $1/n$ and they must be adjusted to make them so in order for them to be useful in practice; however, as $n$ goes to infinity, the proportion $n_m/n$ of observations taken at $x_m$ converges to $\omega_m$. }
Some classical references concerning optimal design theory are \cite{Fedo:Theo:1972} and \cite{Paz:86}.

Since the design (\ref{eq:opt_dis}) depends on the unknown parameter $\boldsymbol{\theta}$ {except} in the case of linear models, it is said {to be} \emph{locally optimal} and can be computed only if a guessed value $\boldsymbol{\theta}_0$ is available. A locally optimal design is usually not robust with respect to different choices of $\boldsymbol{\theta}_0$. To {protect against poor choices of $\boldsymbol{\theta}_0$, one can use} a two stage adaptive procedure where in the first stage $n_1$ observations are recruited according to some design and in the second phase additional $n_2$ data are observed according to a locally optimal design {in which $\boldsymbol{\theta}$ is estimated}  from the first stage data. The whole vector of observations (first and second stage data) are then used to estimate $\boldsymbol{\theta}$ through the maximum likelihood method. The two-stage adaptive design is explained in detail in Section \ref{Sect:twostage}.

The properties of {a multivariate} MLE are studied in Section \ref{sec:AsyProp} assuming model \eqref{eq:Model0} 
in the case that
only the second stage sample size goes to infinity; $n_1$ is assumed to be finite and small. In many different contexts it is quite common to develop a preliminary small pilot study in order to have an idea about the phenomenon under study and then to perform a larger and well developed study on the same subject. Thus, it is practical to assume that $n_1$ is fixed and small, and then asymptotic approximation in the first stage is not adequate.

\section{Two-stage adaptive design and corresponding model}\label{Sect:twostage}
Assume 
that in the first stage a finite number of independent observations, say $n_1<+\infty$, are taken according to a design 
 $$\xi_1=\left\{\!\!\!\begin{array}{ccc}
x_{11} &\cdots & x_{1M_1} \\
\omega_{11} & \cdots & \omega_{1M_1}
\end{array}
\!\!\!  \right\},$$
i.e. $n_{1m}=n_1 \omega_{1m}$ observations are taken at the experimental point $x_{1m}$, for $m=1,\ldots, M_1$.

Let $\{y_{1mj}\}_{1,1}^{M_1,n_{1m}}$  be the first stage observations. An estimate for $\boldsymbol{\theta}$ can be computed maximizing the likelihood corresponding to these first stage observations; the MLE $\boldsymbol{\hat\theta_{n_1}}$ depends on the first stage data  through the complete sufficient statistic $\mathbf{\bar{y}}_{1}=(\bar y_{11},\ldots,\bar y_{1M_1})^T$, where
$\bar y_{1m}= \sum_{j=1}^{n_{1m}} y_{1mj}/n_{1m}$, $m=1,\ldots,M_1$; thus,  $\boldsymbol{\hat\theta_{n_1}}=\boldsymbol{\hat\theta_{n_1}}(\mathbf{\bar{y}}_{1})$.

 In the second stage, $n_2$ independent observations are accrued according to the following local optimum design
\begin{equation}\label{eq:xi2star}
 \xi^*_2=\xi_2^*(\boldsymbol{\hat\theta_{n_1}})=\left\{\!\!\!\begin{array}{ccc}
 x_{21} &\cdots & x_{2M_2} \\
 \omega_{21} & \cdots & \omega_{2M_2}
 \end{array}
 \!\!\!  \right\};
\end{equation}
$\{y_{2mj}\}_{1,1}^{M_2,n_{2m}}$ denotes the second stage observations, where  $n_{2m}$ is obtained by rounding $n_2 \omega_{2m}$ to an integer under the constraint $\sum_{m=1}^{M_2}n_{2m}=n_2$, for $m=1,\ldots, M_2$.

Note that $\xi^*_2$ is a random probability distribution (discrete and finite) since it depends on the first stage observation through $\mathbf{\bar{y}}_{1}$ as $\boldsymbol{\hat\theta_{n_1}}=\boldsymbol{\hat\theta_{n_1}}(\mathbf{\bar{y}}_{1})$; thus, given $\mathbf{\bar{y}}_{1}$, the second stage design  $\xi^*_2$ is determined and  $\{y_{2mj}\}_{1,1}^{M_2,n_{2m}}$ are  $n_2$ conditionally independent observations. In addition,
 it is natural to assume that  second stage observations depend on the first stage information only through $\xi^*_2$. As a consequence,   the observations $\left\{y_{imj}\right\}_{1,1,1}^{2,M_i,n_{im}}$ follow the model
\begin{align}\label{eq:Model}
y_{imj} = \eta(x_{im},\boldsymbol{\theta}) + \varepsilon_{imj},   
\end{align}
where, given $\xi^*_2$, $\{y_{2mj}\}_{1,1}^{M_2,n_{2m}}$ are conditionally independent of  $\{y_{1mj}\}_{1,1}^{M_1,n_{1m}}$, and $\varepsilon_{imj}$ are i.i.d. $\mathscr{N}\left( 0,\sigma^2 \right)$ for any $i,m,j$.


\subsection{Likelihood and the Fisher information matrix}
The likelihood for model (\ref{eq:Model}) is
\begin{align}\label{eq:tlikelihood}
\mathscr{L}_{n}\left(\boldsymbol{\theta}|\mathbf{\bar{y}}_{1},\mathbf{\bar{y}}_{2}, \mathbf{x}_1, \mathbf{x}_2\right) \propto
{\cal L}_{1n_1}\left(\boldsymbol{\theta}|\mathbf{\bar{y}}_{1}, \mathbf{x}_1\right)\cdot {\cal L}_{2n_2}\left(\boldsymbol{\theta}|\mathbf{\bar{y}}_{2}, \mathbf{x}_2\right),
\end{align}
where {$n=n_1+n_2$ and}
$$ {\cal L}_{in_i}\left(\boldsymbol{\theta}|\mathbf{\bar{y}}_{i}, \mathbf{x}_i\right) \propto \exp\left\{ -\frac{1}{2\sigma^{2}}\sum_{m=1}^{M_i} n_{im}\left[ \bar{y}_{im} - \eta(x_{im},\boldsymbol{\theta})\right]^{2}
 \right\},\quad i=1,2;$$
\noindent   $\mathbf{\bar{y}}_{i}=(\bar y_{i1},\ldots,\bar y_{iM_i})^T$ { and }
$\bar y_{im}= n_{im}^{-1} \sum_{j=1}^{n_{im}} y_{imj}$ is the stage $i$ sample mean at the $m$-th dose for $m=1,\ldots,M_i$.

The total score function is
\begin{equation}
\label{eq:totalS}
\mathbf{S}_n = \nabla\ln \mathscr{L}_{n}\left(\boldsymbol{\theta}|\mathbf{\bar{y}}_{1},\mathbf{\bar{y}}_{2}, \mathbf{x}_1, \mathbf{x}_2\right)=
 \mathbf{S}_{1 n_1} + \mathbf{S}_{2 n_2},
 \end{equation}
 where
\begin{align*}
 \mathbf{S}_{i n_i} =\nabla \ln {\cal L}_{in_i}\left(\boldsymbol{\theta}|\mathbf{\bar{y}}_{i}, \mathbf{x}_i\right)= \frac{1}{\sigma^{2}}
\sum_{m=1}^{M_i} n_{im}\left[ \bar{y}_{im} - \eta(x_{im},\boldsymbol{\theta})\right]\nabla  \eta(x_{im},\boldsymbol{\theta})
\end{align*}
represents the score function for the $i$-th stage.

As outlined before, $\mathbf{\bar{y}}_2$ depends on $\mathbf{\bar{y}}_1$   only through $\xi_2^*$ and, given $\mathbf{\bar{y}}_1$, the second stage design $\xi_2^*$ is completely determined. As a consequence, $\textrm{E}_{\mathbf{\bar{y}}_{2}|\mathbf{\bar{y}}_{1}}[\mathbf{S_2}]=0$  and
Fisher information matrix is
\begin{equation}
\textrm{Cov}_{\mathbf{\bar{y}}_{1},\mathbf{\bar{y}}_{2}}[\mathbf{S}_n,\mathbf{S}_n] =
\mathbf{\textrm{E}}_{\mathbf{\bar{y}}_{1}}
[\mathbf{S}_{1n_1}\mathbf{S}_{1n_1}^{\textrm{\sc{T}}}]+
\mathbf{\textrm{E}}_{\mathbf{\bar{y}}_{1},\mathbf{\bar{y}}_{2}}
[\mathbf{S}_{2n_2}\mathbf{S}_{2n_2}^{\textrm{\sc{T}}}], 
\label{Fisher}
\end{equation}
where
\begin{align*}
\mathbf{\textrm{E}}_{\mathbf{\bar{y}}_{1}}
[\mathbf{S}_{1n_1}\mathbf{S}_{1n_1}^{\textrm{\sc{T}}}]
&=
\frac{1}{\sigma^{2}}
\sum_{m=1}^{M_1} n_{1m}
\nabla \eta(x_{1m},\boldsymbol{\theta})
\nabla \eta(x_{1m},\boldsymbol{\theta})^{\textrm{\sc T}};\\
\mathbf{\textrm{E}}_{\mathbf{\bar{y}}_{1},\mathbf{\bar{y}}_{2}}
[\mathbf{S}_{2n_2}\mathbf{S}_{2n_2}^{\rm T}]
&=
\mathbf{\textrm{E}}_{\mathbf{\bar{y}}_{1}}
\mathbf{\textrm{E}}_{\mathbf{\bar{y}}_{2}|\mathbf{\bar{y}}_{1}}
[\mathbf{S}_{2n_2}\mathbf{S}_{2n_2}^{\textrm{\sc{T}}}];\\
\mathbf{\textrm{E}}_{\mathbf{\bar{y}}_{2}|\mathbf{\bar{y}}_{1}}
[\mathbf{S}_{2n_2}\mathbf{S}_{2n_2}^{\textrm{\sc{T}}}]
&=
\frac{1}{\sigma^{2}}
\sum_{m=1}^{M_2} n_{2m}
\nabla \eta(x_{2m},\boldsymbol{\theta})
\nabla \eta(x_{2m},\boldsymbol{\theta})^{\textrm{\sc T}}.
\end{align*}
Now, the per-subject information can be written as
\begin{align*}
 \frac{1}{n}\mathbf{\textrm{Cov}}_{\mathbf{\bar{y}}_{1},\mathbf{\bar{y}}_{2}}[\mathbf{S_n}\mathbf{S}_n^{\textrm{\sc{T}}}]
& = \frac{1}{n\sigma^{2}}
 \left\{
 \sum_{m=1}^{M_1} n_{1m} \nabla \eta(x_{1m},\boldsymbol{\theta})
\nabla \eta(x_{1m},\boldsymbol{\theta})^{\textrm{\sc{T}}}
\right.\\
&\hspace{2cm} +\left.
 \mbox{E}_{\mathbf{\bar{y}}_{1}}\left[\sum_{m=1}^{M_2}   n_{2m}
\nabla \eta(x_{2m},\boldsymbol{\theta})
\nabla \eta(x_{2m},\boldsymbol{\theta})^{\textrm{\sc{T}}}
\right]
\right\},
\end{align*}
where $M_2$, ${x}_{2m}$ and ${n}_{2m}$ are  random variables,  defined by the onto transformation \eqref{eq:xi2star} of $\mathbf{\bar{y}}_{1}$.

{N}ote that, as $n_2\rightarrow \infty$ (and thus $n\rightarrow \infty$), the per-subject information converges almost surely to 
 \begin{equation}
  \frac{1}{\sigma^2}\mbox{E}_{\mathbf{\bar{y}}_{1}}\left[\int_{\cal X}   
 \nabla \eta(x,\boldsymbol{\theta})
 \nabla \eta(x,\boldsymbol{\theta})^{\textrm{\sc{T}}}\,d\xi_2^*(x)
 \right].
 \label{per-subjectInf}
 \end{equation}


\section{Asymptotic Properties}\label{sec:AsyProp}
One needs an approximation to the asymptotic distribution of the final MLE $\widehat{\boldsymbol{\theta}}_{n}$ {that} may be used for inference at the end of the study, where $n=n_1+n_2$ is the total number of observations.
The classical approach is to assume that both $n_{1}$ and $n_{2}$ are large (see for instance \cite{Pronz} and \cite{Dette:2012}).  {For a closer approximation to  many experimental situations, } assume here a fixed first stage sample size $n_{1}$ and a large second stage sample size $n_{2}$.

{N}ote that if the experimental conditions in model \eqref{eq:Model0} are taken according to an experimental design $\xi$, then, {by} the law of large numbers,
\begin{equation}\label{LLN}
\dfrac{1}{n}\sum_{i=1}^{n}\eta(x_{i},\boldsymbol{\theta})\,\eta(x_{i},\boldsymbol{\theta}_1)\overset{P}\longrightarrow \int \eta(x,\boldsymbol{\theta})\,\eta(x,\boldsymbol{\theta}_1) d\xi.
\end{equation}
 In order to prove the consistency of $\widehat{\boldsymbol{\theta}}_{n}$  assume the following:

\begin{assumption}\label{A1}
The model is identifiable: if ${\boldsymbol{\theta}}_1\neq {\boldsymbol{\theta}}_2$, then $\eta(x,\boldsymbol{\theta}_1)\neq\eta(x,\boldsymbol{\theta}_2)$.
\end{assumption}

\begin{assumption}\label{A2}
The convergence \eqref{LLN} is uniform for all $\boldsymbol{\theta},\boldsymbol{\theta}_1 \in \boldsymbol{\Theta}$,
that is, for any $\delta>0$,
\begin{equation*}\label{ULLN}
P\left(\sup_{\boldsymbol{\theta},\boldsymbol{\theta}_1 \in \boldsymbol{\Theta}}\left|\dfrac{1}{n}\sum_{i=1}^{n}\eta(x_{i},\boldsymbol{\theta})\,\eta(x_{i},\boldsymbol{\theta}_1)- \int \eta(x,\boldsymbol{\theta})\,\eta(x,\boldsymbol{\theta}_1) d\xi\right|>\delta\right)\longrightarrow 0.\\[12pt]
\end{equation*}
\end{assumption}

\begin{theorem}
	\label{Th:1}
	Let $\widehat{\boldsymbol{\theta}}_{n}$ be the MLE maximing the total likelihood  \eqref{eq:tlikelihood}. Then $$\widehat{\boldsymbol{\theta}}_{n}\overset{P}\longrightarrow {\boldsymbol{\theta}}^t,$$
	where ${\boldsymbol{\theta}}^t$ denotes the true unknown value of ${\boldsymbol{\theta}}$.
\end{theorem}

\noindent \textbf{Proof.} Observe that $\widehat{\boldsymbol{\theta}}_{n}$ maximizes  \eqref{eq:tlikelihood} if and only if it minimizes the average square{d} errors
\begin{eqnarray}\label{SS}
{\cal A}_n(\boldsymbol{\theta})&=&\dfrac{1}{n}\sum_{i=1}^{n_1}[y_{i1}-\eta(x_{i1},\boldsymbol{\theta})]^2+\dfrac{1}{n}\sum_{i=1}^{n_2}[y_{i2}-\eta(x_{i2},\boldsymbol{\theta})]^2.
\end{eqnarray}
{To} prove that 
\begin{equation}\label{ASEn}
\sup_{\boldsymbol{\theta} \in \boldsymbol{\Theta}} |{\cal A}_n(\boldsymbol{\theta})-{\cal A}(\boldsymbol{\theta})| \overset{P}\longrightarrow 0,
\end{equation}
where
\begin{equation}\label{ASE}
{\cal A}(\boldsymbol{\theta})=\sigma^2+\int [\eta(x,\boldsymbol{\theta}^t)-\eta(x,\boldsymbol{\theta})]^2 d\xi^*_2,
\end{equation}
Rewrite {\eqref{SS} as}
\begin{eqnarray}
{\cal A}_n(\boldsymbol{\theta})&=& \dfrac{1}{n}\sum_{i=1}^{n_1}[y_{i1}-\eta(x_{i1},\boldsymbol{\theta})]^2
+\dfrac{1}{n}\sum_{i=1}^{n_2}[y_{i2}-\eta(x_{i2},\boldsymbol{\theta}^t)+\eta(x_{i2},\boldsymbol{\theta}^t)-\eta(x_{i2},\boldsymbol{\theta})]^2 \nonumber\\
&=& A_n(\boldsymbol{\theta})+B_n(\boldsymbol{\theta}^t)+C_n(\boldsymbol{\theta})+D_n(\boldsymbol{\theta}),\nonumber
\end{eqnarray}
where 
\begin{align*}
A_n(\boldsymbol{\theta})&=\dfrac{1}{n}\sum_{i=1}^{n_1}[y_{i1}-\eta(x_{i1},\boldsymbol{\theta})]^2;\\
B_n(\boldsymbol{\theta}^t)&=\dfrac{1}{n}\sum_{i=1}^{n_2}[y_{i2}-\eta(x_{i2},\boldsymbol{\theta}^t)]^2; \\
C_n(\boldsymbol{\theta})&=\dfrac{2}{n}\sum_{i=1}^{n_2}[y_{i2}-\eta(x_{i2},\boldsymbol{\theta}^t)][\eta(x_{i2},\boldsymbol{\theta}^t)-\eta(x_{i2},\boldsymbol{\theta})]; \\
D_n(\boldsymbol{\theta})&=\dfrac{1}{n}\sum_{i=1}^{n_2}[\eta(x_{i2},\boldsymbol{\theta}^t)-\eta(x_{i2},\boldsymbol{\theta})]^2.
\end{align*}
{It follows that}
\begin{enumerate}
	\item $\sup_{\boldsymbol{\theta} \in \boldsymbol{\Theta}} |A_n(\boldsymbol{\theta})| \overset{P}\longrightarrow 0$ because $n_1$ is finite;
	\item $B_n(\boldsymbol{\theta}^t)=\dfrac{1}{n}\sum_{i=1}^{n_2}\varepsilon_{i2}^2\overset{P}\longrightarrow \sigma^2$ because the $\{\varepsilon_{i2}\}_{i=1}^\infty$ is a sequence of i.i.d. random variables $\sim{\cal N}(0;\sigma^2)$;
	\item {T}he random variables
	\begin{equation*}
	[y_{i2}-\eta(x_{i2},\boldsymbol{\theta}^t)][\eta(x_{i2},\boldsymbol{\theta}^t)-\eta(x_{i2},\boldsymbol{\theta})] = \epsilon_{i2}[\eta(x_{i2},\boldsymbol{\theta}^t)-\eta(x_{i2},\boldsymbol{\theta})]
	\end{equation*}
	are i.i.d. conditionally {on} $\xi^*_2$ and 
	\begin{equation*}
	E[\epsilon_{i2}(\eta(x_{i2},\boldsymbol{\theta}^t)-\eta(x_{i2},\boldsymbol{\theta}))|\xi^*_2]=(\eta(x_{i2},\boldsymbol{\theta}^t)-\eta(x_{i2},\boldsymbol{\theta}))E[\epsilon_{i2}]=0.
	\end{equation*}
	Hence, from the conditional law of large numbers (see, for instance, \cite[Theorem 7]{PRao}) and {because} $\eta$ is continuous on the compact set $\Theta$,
	$$\sup_{\boldsymbol{\theta} \in \boldsymbol{\Theta}} |C_n(\boldsymbol{\theta})| \overset{P}\longrightarrow 0.$$	
	\item Notice that $$D_n(\boldsymbol{\theta})=\dfrac{1}{n}\sum_{i=1}^{n_2}\eta(x_{i2},\boldsymbol{\theta}^t)^2+\dfrac{1}{n}\sum_{i=1}^{n_2}\eta(x_{i2},\boldsymbol{\theta})^2-\dfrac{2}{n}\sum_{i=1}^{n_2}\eta(x_{i2},\boldsymbol{\theta}^t)\eta(x_{i2},\boldsymbol{\theta});$$
	hence, {from} the conditional law of large numbers and  Assumption \ref{A2}, 
	\begin{equation*}
	P(\sup_{\boldsymbol{\theta} \in \boldsymbol{\Theta}}|D_n(\boldsymbol{\theta})-D(\boldsymbol{\theta}) |>\delta|\xi^*_2)\longrightarrow 0,
	\end{equation*}
	a.s. for any $\delta>0$, where
	\begin{align*}
	D(\boldsymbol{\theta})&= \int\eta(x,\boldsymbol{\theta}^t)^2 d\xi^*_2+\int\eta(x_{i2},\boldsymbol{\theta})^2 d\xi^*_2-2\int \eta(x_{i2},\boldsymbol{\theta}^t)\eta(x_{i2},\boldsymbol{\theta}) d\xi^*_2\\
	&=\int [\eta(x,\boldsymbol{\theta}^t)-\eta(x,\boldsymbol{\theta})]^2 d\xi^*_2.
	\end{align*}
	It follows that
	\begin{equation*}
	E[P(\sup_{\boldsymbol{\theta} \in \boldsymbol{\Theta}}|D_n(\boldsymbol{\theta})-D(\boldsymbol{\theta}) |>\delta|\xi^*_2)]
	=P(\sup_{\boldsymbol{\theta} \in \boldsymbol{\Theta}}|D_n(\boldsymbol{\theta})-D(\boldsymbol{\theta}) |>\delta)\longrightarrow 0.
	\end{equation*}
\end{enumerate}
{Statement} \eqref{ASEn} {is} proved. $\boldsymbol{\theta}^t$ is the unique minimum of ${\cal A}(\boldsymbol{\theta})$   as a consequence of Assumption \ref{A1} and hence the thesis follows.

\begin{theorem} \label{thm:UQ} For model (\ref{eq:Model}) with $\xi_{2}^*$ defined in \eqref{eq:xi2star}  and $M(\cdot,\cdot)$ defined in \eqref{M},
\begin{equation}\label{eq:FixedStageOne}
\sqrt{n}\left(\widehat{\boldsymbol{\theta}}_{n} - \boldsymbol{\theta}^{t}\right) \xrightarrow{\mathscr{D}}
\sigma\; M(\xi^*_2,\boldsymbol{\theta}^t)^{-1/2}\;\boldsymbol{Z}
\end{equation}
as $n_{2} \rightarrow \infty$,
where $\boldsymbol{Z}$ is a $(p+1)$-dimensional standard normal random vector independent of the random matrix $M(\xi^*_2,\boldsymbol{\theta}^t)$.
\end{theorem}

\noindent \textbf{Proof.} Let $S_{n}^j$ be the $j$-the component of the total score function $\mathbf{S}_n$  in \eqref{eq:totalS}. From the expansion of  $\mathbf{S}_n(\boldsymbol{\theta})$ around the true value ${\boldsymbol{\theta}}^t$ 
we obtain, for any {parameter} $j=0, ..., p$,
\begin{align*}
S_{n}^j(\widehat{\boldsymbol{\theta}}_{n})=
{S_{1n_1}^j}({\boldsymbol{\theta}}^t)+
{S_{2n_2}^j}({\boldsymbol{\theta}}^t)&+
\sum\limits_{k=0}^{p}(\widehat\theta_{nk}-\theta^t_k) \dot{S}_n^{jk}({\boldsymbol{\theta}}^t)\\ & +
\frac 1 2 \sum_{k=0}^{p}\sum_{l=0}^{p}(\widehat\theta_{nk}-\theta^t_k)
(\widehat\theta_{nl}-\theta^t_l) \ddot{S}_n^{jkl}({\boldsymbol{\theta}}^*), \end{align*}
where
\begin{align*}\dot{S}_n^{jk}({\boldsymbol{\theta}})&=
\frac{\partial^2}{\partial\theta_j\partial\theta_k}\ln{\cal L}_{1n_1}(\boldsymbol{\theta})+
\frac{\partial^2}{\partial\theta_j\partial\theta_k}\ln{\cal L}_{2n_2}(\boldsymbol{\theta})=
\dot{S}_{1n_1}^{jk}({\boldsymbol{\theta}})+\dot{S}_{2n_2}^{jk}({\boldsymbol{\theta}}),\\
\ddot{S}_n^{jkl}({\boldsymbol{\theta}})&=
\frac{\partial^3}{\partial\theta_j\partial\theta_k\partial\theta_l}\ln{\cal L}_{1n_1}(\boldsymbol{\theta})+
\frac{\partial^3}{\partial\theta_j\partial\theta_k\partial\theta_l}\ln{\cal L}_{2n_2}(\boldsymbol{\theta})=
\ddot{S}_{1n_1}^{jkl}({\boldsymbol{\theta}})+\ddot{S}_{2n_2}^{jkl}({\boldsymbol{\theta}})\end{align*}
and ${\boldsymbol{\theta}}^*$ is a point between $\widehat{\boldsymbol{\theta}}_{n}$ and ${\boldsymbol{\theta}}^t$. 
Since ${S_n^j}(\widehat{\boldsymbol{\theta}}_{n})=0$, 
\begin{align*}
\dfrac 1 {\sqrt n}{S_{2n_2}^j}({\boldsymbol{\theta}}^t) = 
&-\dfrac 1 {\sqrt n}   \left[{S_{1n_1}^j}({\boldsymbol{\theta}}^t)+
\sum_{k=0}^{p} (\widehat\theta_{nk}-\theta^t_k) \dot{S}_{1n_1}^{jk}({\boldsymbol{\theta}}^t)\right]\\
&-
\dfrac 1 {\sqrt n} \left[\sum_{k=0}^{p}\sum_{l=0}^{p}(\widehat\theta_{nk}-\theta^t_k)
(\widehat\theta_{nl}-\theta^t_l) \ddot{S}_{1n_1}^{jkl}({\boldsymbol{\theta}}^*)\right] \\
&+
\sqrt n \sum_{k=0}^{p}(\widehat\theta_{nk}-\theta^t_k)\left[-\frac 1 n \dot{S}_{2n_2}^{jk}({\boldsymbol{\theta}}^t)- \dfrac 1 {2n}
\sum_{l=0}^{p}
(\widehat\theta_{nl}-\theta^t_l) \ddot{S}_{2n_2}^{jkl}({\boldsymbol{\theta}}^*) \right].
\end{align*}
From the consistency proved in Theorem \ref{Th:1}, 
 $S_{2n_2}^j({\boldsymbol{\theta}}^t)/ \sqrt n$
is asymptotically equivalent to 
$$
\sqrt n \sum_{k=0}^{p}(\widehat\theta_{nk}-\theta^t_k)\left[-\frac 1 n \dot{S}_{2n_2}^{jk}({\boldsymbol{\theta}}^t)- \dfrac 1 {2n}
\sum_{l=0}^{p}
(\widehat\theta_{nl}-\theta^t_l) \ddot{S}_{2n_2}^{jkl}({\boldsymbol{\theta}}^*) \right],
$$
in the sense that their difference converges in probability to zero. In matrix notation, let
$$
\mathbf{\dot{S}}_{2n_2}({\boldsymbol{\theta}}^t)=\left\{ \dot{S}_{2n_2}^{jk}({\boldsymbol{\theta}}^t)\right\}_{(jk)} 
\quad {\rm and} \quad
\mathbf{\ddot{S}}_{2n_2}^{(l)}({\boldsymbol{\theta}}^*)=\left\{ \ddot{S}_{2n_2}^{jkl}({\boldsymbol{\theta}}^*)\right\}_{(jk)}, \quad j,k=0,\ldots,p,
$$ 
then
\begin{equation}\label{eq:as_eq}
\dfrac 1 {\sqrt n}\mathbf{S}_{2n_2}({\boldsymbol{\theta}}^t)
\quad \mbox{and}\quad
\left[-\frac{1}{n} \mathbf{\dot{S}}_{2n_2}({\boldsymbol{\theta}}^t) -\frac{1}{2n}\sum_{l=0}^{p} (\widehat\theta_{nl}-\theta^t_l)  \mathbf{\ddot{S}}_{2n_2}^{(l)}({\boldsymbol{\theta}}^*) \right] \sqrt{n}\left(\widehat{\boldsymbol{\theta}}_{n} - \boldsymbol{\theta}^{t}\right),
 \end{equation}
are asymptotically equivalent. 

Now,
\begin{eqnarray}\label{primopezzo}
\nonumber \dfrac 1 {\sqrt n}\mathbf{S}_{2n_2}({\boldsymbol{\theta}}^t)&=& \frac{1}{\sigma^{2}}\dfrac 1 {\sqrt n} \, \sum_{i=1}^{n_2}\, [y_{2i}-\eta(x_{2i},\boldsymbol{\theta}^t)]\nabla\eta(x_{2i},\boldsymbol{\theta}^t)\\
 &=& \frac{1}{\sigma} \dfrac 1 {\sqrt n} \, \sum_{i=1}^{n_2}\,\frac{1}{\sigma} \varepsilon_{2i}\;\nabla\eta(x_{2i},\boldsymbol{\theta}^t) 
 \end{eqnarray}
is a zero-mean, square integrable, {\it martingale difference array} with respect to the filtration ${\cal F}_{0}=\sigma(\mathbf{\bar{y}}_{1})$, ${\cal F}_{1}=\sigma(\mathbf{\bar{y}}_{1}, \varepsilon_{21})$, ..., ${\cal F}_{n_2}=\sigma(\mathbf{\bar{y}}_{1}, \varepsilon_{21}, \ldots, \varepsilon_{2n_2})$,
according to the definition in \cite{HallHeyde}.

 It follows from \cite[Theorem 3.2]{HallHeyde} that
\begin{equation}\label{array}
\dfrac 1 {\sqrt n}\mathbf{S}_{2n_2}({\boldsymbol{\theta}}^t) \xrightarrow{\mathscr{D}} \frac{1}{\sigma} M(\xi^*_2,\boldsymbol{\theta}^t)^{1/2}\boldsymbol{Z} \quad \mbox{(stably)}
\end{equation}
as $n_{2} \rightarrow \infty$,
where $\boldsymbol{Z}$ is a $(p+1)$-dimensional standard normal random vector independent of the random matrix $M(\xi^*_2,\boldsymbol{\theta}^t)$. 
Note that Assumptions 3.18 and 3.20 of \cite[Theorem 3.2]{HallHeyde} are easily verified, while {Assumption}  3.19 becomes
\begin{equation}\label{eq:3.19}
\frac{1}{\sigma^{4}}\dfrac 1 {n} \sum_{i=1}^n  \varepsilon_{2i}^2\;\nabla\eta(x_{2i},\boldsymbol{\theta}^t) \nabla \eta(x_{2m},\boldsymbol{\theta}^t)^T \overset{P}\longrightarrow \frac{1}{\sigma^2} M(\xi^*_2,\boldsymbol{\theta}^t).
\end{equation}
To obtain the \eqref{eq:3.19}, the conditional law of large numbers \cite[Theorem 7]{PRao} can be applied: conditional on $\sigma(\mathbf{\bar{y}}_{1})$,
\begin{eqnarray}\label{etaq}
 \dfrac 1 {n} \sum_{i=1}^n  \varepsilon_{2i}^2\;\nabla\eta(x_{2i},\boldsymbol{\theta}^t) \nabla \eta(x_{2i},\boldsymbol{\theta}^t)^T &\overset{P}\rightarrow&  E[\varepsilon_{2}^2\;\nabla\eta(x_{2},\boldsymbol{\theta}^t) \nabla \eta(x_{2},\boldsymbol{\theta}^t)^T|\mathbf{\bar{y}}_{1}]\nonumber \\
 \nonumber
&=&  E[\varepsilon_{2}^2|\mathbf{\bar{y}}_{1}]\cdot E[\nabla\eta(x_{2},\boldsymbol{\theta}^t) \nabla \eta(x_{2},\boldsymbol{\theta}^t)^T|\mathbf{\bar{y}}_{1}]\\ 
&=& \sigma^2 \int_{\cal X} \nabla \eta(x,\boldsymbol{\theta}^t)
\nabla \eta(x,\boldsymbol{\theta}^t)^T d\xi^*_2(x);
\end{eqnarray}
averaging on the conditional probability, the convergence \eqref{etaq} mantains also unconditionally. 

 As a consequence of \eqref{array}, as shown in \cite[(vi) in \S 3.2]{HallHeyde}, since $M(\xi^*_2,\boldsymbol{\theta}^t)$ is ${\cal F}_{n_2}$-measurable for all $n_2$,
 \begin{equation}
 \sigma\, M(\xi^*_2,\boldsymbol{\theta}^t)^{-1/2}\dfrac 1 {\sqrt n}\mathbf{S}_{2n_2}({\boldsymbol{\theta}}^t) \xrightarrow{\mathscr{D}} \boldsymbol{Z},
 \end{equation}
 where $\boldsymbol{Z}$ is a $(p+1)$-dimensional standard normal random vector independent of the random matrix $M(\xi^*_2,\boldsymbol{\theta}^t)$.
 Thus, \eqref{eq:as_eq} {provides that} also
\begin{equation}
Q_n=\sigma\, M(\xi^*_2,\boldsymbol{\theta}^t)^{-1/2}\left[-\frac{1}{n} \mathbf{\dot{S}}_{2n_2}({\boldsymbol{\theta}}^t) -\frac{1}{2n}\sum_{l=0}^{p} (\widehat\theta_{nl}-\theta^t_l)  \mathbf{\ddot{S}}_{2n_2}^{(l)}({\boldsymbol{\theta}}^*) \right] 
\sqrt{n}\left(\widehat{\boldsymbol{\theta}}_{n} - \boldsymbol{\theta}^{t}\right) 
 \xrightarrow{\mathscr{D}} \boldsymbol{Z},
\end{equation}
from Slutsky's theorem.
Moreover,
\begin{align}\label{secondopezzo}
 -\frac 1 n \dot{\mathbf S}_{2n_2}({\boldsymbol{\theta}}^t) \xrightarrow{\mathscr{P}}
 \frac{1}{\sigma^2}\sum_{m=1}^{M_2} \omega_{2m} \nabla \eta(x_{2m},\boldsymbol{\theta}^t)\nabla \eta(x_{2m},\boldsymbol{\theta}^t)^T
=\frac{1}{\sigma^2} M(\xi^*_2,\boldsymbol{\theta}^t)
\end{align}
 because the $jk$-th element of the matrix $- \dot{\mathbf S}_{2n_2}({\boldsymbol{\theta}}^t)/n$ satisfies
\begin{align}\label{eq:utile}
- &\dfrac 1 n \dot{S}_{2n_2}^{jk}({\boldsymbol{\theta}}^t) \nonumber \\
&= \frac{1}{\sigma^2}\sum_{m=1}^{M_2} \dfrac{n_{2m}} { n} \left[\dfrac{\partial \eta(x_{2m},\boldsymbol{\theta}^t)}{\partial \theta_j}\cdot\dfrac{\partial \eta(x_{2m},\boldsymbol{\theta}^t)}{\partial \theta_k}- \dfrac{\partial^2 \eta(x_{2m},\boldsymbol{\theta}^t)}{\partial \theta_j \partial \theta_k} [\bar{y}_{2m} - \eta(x_{2m},\boldsymbol{\theta}^t)] \right]
\end{align}
and the {last} right term of equation \eqref{eq:utile} converges in probability to zero by the conditional law of large numbers. \\
Now,
\begin{equation}
\sqrt{n}\left(\widehat{\boldsymbol{\theta}}_{n} - \boldsymbol{\theta}^{t}\right) =
R_n  \cdot  Q_n,
\end{equation}
where 
\begin{equation*}
R_n:=\left[-\frac{1}{n} \mathbf{\dot{S}}_{2n_2}({\boldsymbol{\theta}}^t) -\frac{1}{2n}\sum_{l=0}^{p} (\widehat\theta_{nl}-\theta^t_l)  \mathbf{\ddot{S}}_{2n_2}^{(l)}({\boldsymbol{\theta}}^*) \right]^{-1}
\cdot \dfrac{1}{\sigma}M(\xi^*_2,\boldsymbol{\theta}^t)^{1/2}.
\end{equation*}
From \eqref{secondopezzo} and from the consistency proved in Theorem \ref{Th:1} (assuming the standard regularity conditions {needed} for  $\frac{1}{n}\ddot{S}_{2n_2}^{jkl}({\boldsymbol{\theta}}^*)$ {to be} bounded in probability), 
$$
R_n  \xrightarrow{\mathscr{P}} \sigma M(\xi^*_2,\boldsymbol{\theta}^t)^{-1/2}.
$$
Since the limits in the distributions of $R_n$ and $Q_n$ are independent, $(R_n,Q_n)$ converges to $[\sigma M(\xi^*_2,\boldsymbol{\theta}^t)^{-1/2},\boldsymbol{Z}]$, and hence,  $R_n\cdot Q_n$ $\xrightarrow{\mathscr{D}}$ $\sigma M(\xi^*_2,\boldsymbol{\theta}^t)^{-1/2} \boldsymbol{Z}$ from Slutsky's theorem, obtaining the thesis.
\bigskip

\begin{corollary}\label{theo:asvar}
The asymptotic variance of $\sqrt{n}\left(\widehat{\boldsymbol{\theta}}_{n} - \boldsymbol{\theta}^{t}\right)$ is
$$
\sigma^2\mbox{E}_{\mathbf{\bar{y}}_{1}}\!\left[ \left(
\int_{\cal X} \nabla \eta(x,\boldsymbol{\theta}^t) \nabla \eta(x,\boldsymbol{\theta}^t)^{\textrm{\sc{T}}}\,d\xi_2^*(x) 
 \right)^{-1}\right] 
$$
\end{corollary}

\noindent \textbf{Proof.}
From (\ref{eq:FixedStageOne}) 
\begin{align}
\mbox{AsVar}\left[\sqrt{n}\left(\widehat{\boldsymbol{\theta}}_{n} - \boldsymbol{\theta}^{t}\right)\right] 
& =
\sigma^2\cdot \mbox{Var} \left[
M(\xi^*_2,\boldsymbol{\theta}^t)^{-1/2}\;\boldsymbol{Z}
\right] \nonumber\\
&=
\sigma^2 \left\{ \mbox{Var}_{\mathbf{\bar{y}}_{1}}  \mbox{E}_{\boldsymbol{Z}}  
 \left[
M(\xi^*_2,\boldsymbol{\theta}^t)^{-1/2}\;\boldsymbol{Z}\;\;\bigg|\mathbf{\bar{y}}_{1}
\right] \right.\nonumber\\
&\qquad  \qquad +
\left. \mbox{E}_{\mathbf{\bar{y}}_{1}}\mbox{Var}_{\boldsymbol{Z}}
\left[
M(\xi^*_2,\boldsymbol{\theta}^t)^{-1/2}\;\boldsymbol{Z}\;\;\bigg|\mathbf{\bar{y}}_{1}
\right]
\right\}.
\label{Q}
\end{align}
Since $\mbox{E}_{Z} ( \boldsymbol{Z}|\mathbf{\bar{y}}_{1})=\mbox{E}_{Z} ( \boldsymbol{Z})= \boldsymbol{0}$, 
the first term in the brackets of (\ref{Q})  vanishes {and in the second term} 
\begin{align*}
\mbox{Var}_{\mathbf{\bar{y}}_{1}}  \mbox{E}_{\boldsymbol{Z}}  
\left[
M(\xi^*_2,\boldsymbol{\theta}^t)^{-1/2}\;\boldsymbol{Z}\;\;\bigg|\mathbf{\bar{y}}_{1}
\right] 
&=
\mbox{Var}_{\mathbf{\bar{y}}_{1}}\!\left[\! M(\xi^*_2,\boldsymbol{\theta}^t)^{-1/2}\cdot 
\mbox{E}_{Z}\!(\boldsymbol{Z})\right] =\boldsymbol{0}.
\end{align*}
Denote by $\boldsymbol{I}$ the identity matrix; taking into account that 
$\mbox{Var}_{Z}\!(\boldsymbol{Z}|\mathbf{\bar{y}}_{1})=\mbox{Var}_{Z}\!(\boldsymbol{Z})=\boldsymbol{I}$, 
the second term in (\ref{Q}) is 
\begin{align*}
\mbox{E}_{\mathbf{\bar{y}}_{1}}\mbox{Var}_{\boldsymbol{Z}}
\left[
M(\xi^*_2,\boldsymbol{\theta}^t)^{-1/2}\;\boldsymbol{Z}\;\;\bigg|\mathbf{\bar{y}}_{1} 
\right]
& = 
\mbox{E}_{\mathbf{\bar{y}}_{1}}\! \left[M(\xi^*_2,\boldsymbol{\theta}^t)^{-1/2}\; \mbox{Var}_Z(\boldsymbol{Z}) \; M(\xi^*_2,\boldsymbol{\theta}^t)^{-1/2}
\right]\\
&=
\mbox{E}_{\mathbf{\bar{y}}_{1}} \left[M(\xi^*_2,\boldsymbol{\theta}^t)^{-1} 
\right],
\end{align*}
and from here the thesis follows.

\bigskip

{\bf Remark.} Compare the asymptotic variance obtained in Corollary \ref{theo:asvar} with the inverse of \eqref{per-subjectInf}, to see that the standard equality between the asymptotic variance of the MLE and the inverse of the per-subject information matrix does not maintain in this context. However,   the asymptotic variance expression obtained in Corollary \ref{theo:asvar}  justifies choosing a design for the second stage by maximizing a concave function of $ M(\xi^*_2,\boldsymbol{\theta}^t)$ as it is commonly done.

\section{Example and simulations: a dose-response model}\label{Sect:simulations}

This section presents a simulation study to compare the two-stage adaptive design with a fixed design in terms of precise estimation.

More specifically, assume 
that a guessed value $\boldsymbol{\theta}_0=(\theta_{0,0}, ..., \theta_{p,0})$ for $\boldsymbol{\theta}$ is available, for instance from  an expert opinion.
Initially we take $n_1$ observations according to a locally D-optimal design
 $$\xi^*_1(\boldsymbol{\theta}_0)=\left\{\!\!\!\begin{array}{ccc}
 x_{11} &\cdots & x_{1M_1} \\
 \omega_{11} & \cdots & \omega_{1M_1}
 \end{array}
 \!\!\!  \right\},$$
 and then:
 \begin{itemize}
 	\item[-] in the fixed design we add $n_2$ observations  according to the same $\xi^*_1(\boldsymbol{\theta}_0)$, independently on the first stage;
 	\item[-] in the adaptive design, instead, 
 	add $n_2$ observations according to the locally optimal design \eqref{eq:xi2star} (with D-optimality).
 \end{itemize}
 In other words, both procedures start with the same fixed optimal design; the fixed continues with this while the adaptive adapts.
 
\subsection{The locally D-optimal design under the $E{max}$ model}
As an example, simulations are performed under the $E{max}$ model, which is well-characterized in the literature and it is  frequently used for dose-response designs in clinical trials, as well as in agriculture and in environmental experiments.  It has  the form \eqref{eq:Model0} with the nonlinear mean function 
\begin{equation}\label{eq:emax}
\eta(x,\boldsymbol{\theta})=\theta_0+\theta_1\,\dfrac{x}{x+\theta_2},
\end{equation}
where $x\in {\cal X}=[a,b]$, $0\le a <b$;  $\theta_0$ represents the response when the dose is zero; $\theta_1$ is the maximum effect attributable to the drug; and $\theta_2$ is the dose which produces the half of the maximum effect.

The locally D-optimal design $\xi^*_D$ for the $Emax$ model is analytically found in \cite{Dette:2010}: 

\begin{equation}
\xi^*_D(\theta_2)=\left\{\!\!\!\begin{array}{ccc}
a & x^*(\theta_2) & b \\
1/3 & 1/3 & 1/3
\end{array}
\!\!\!  \right\},
\label{D-opt}
\end{equation}
where the interior support point $x^*(\theta_2)$ is
\begin{equation}
x^*(\theta_2)=\dfrac{b(a+\theta_2)+a(b+\theta_2)}{(a+\theta_2)+(b+\theta_2)}.
\label{middle_point}
\end{equation}


\bigskip



\subsection{Simulations of MLEs efficiencies}
\begin{figure}
	\centering
	\subfigure[$\theta_0^t=25,50$; $\sigma=0.1$]{
		\includegraphics[width=0.49\textwidth]{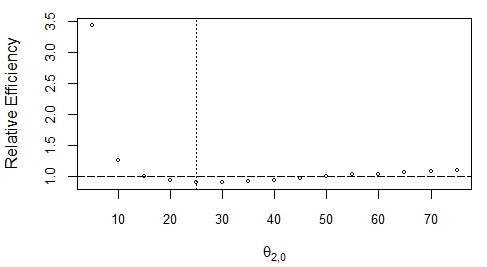}
		\includegraphics[width=0.49\textwidth]{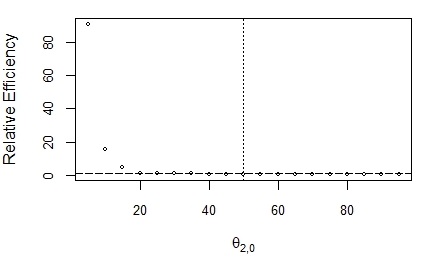}	
	}
	\subfigure[$\theta_0^t=25,50$; $\sigma=0.25$]{
		\includegraphics[width=0.49\textwidth]{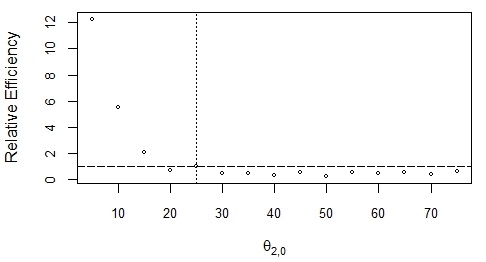}
		\includegraphics[width=0.49\textwidth]{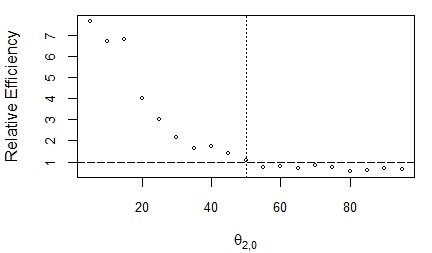}	
	}
	\subfigure[$\theta_0^t=25,50$; $\sigma=0.5$]{
		\includegraphics[width=0.49\textwidth]{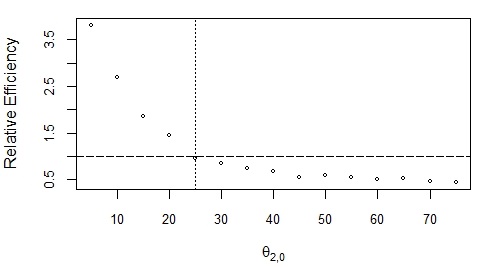}
		\includegraphics[width=0.49\textwidth]{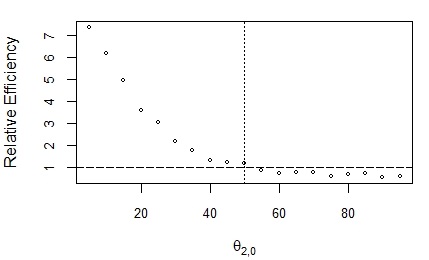}	
	}
	\caption{
	Relative efficiency versus $\theta_{2,0}$ under the Emax model. Relative efficiency  is the MSE of $\widehat{\boldsymbol{\theta}}_{n}$ under the adaptive procedure divided by the MSE under the fixed procedure. The vertical line represents the value of $\theta^t_0: \theta^t_0= 25$ on the left and $\theta^t_0: \theta^t_0= 50$ on the right.}
	\label{relEff}
\end{figure}

  {To}  compare the precision of  the {MLEs} $\widehat{\boldsymbol{\theta}}_{n}$ obtained from the fixed and from the adaptive procedures, we compute the corresponding MSEs and their relative efficiency. In these simulations  $n_1=27$ and  $n_2=270$.  The results are obtained with the  $R$ package developed by \cite{Bornkamp2010} and are based on 10.000 repetitions.  The domain of the non-linear parameter $\theta_2$ is $[0.015;1500]$ to ensure the existence of the MLE; since $\xi^*_D$ does not depends on $\theta_0$ and $\theta_1$, the true values of the linear parameters are fixed at $\theta_0^t=2$ and $\theta_1^t=0.467$, as in \cite{Dette:2012}. 

  Figure 1 display the relative efficiency of the adaptive design with respect to the fixed, for different true values $\theta_2^t=25$ and $\theta_2^t= 50$ and for different nominal values $\theta_{2,0}$, varying on the $x$-axis around the true value. The standard deviation is assumed to be $\sigma=0.1, 0.25, 0.5$. To give an idea of the order of magnitude of the MSE, some values are reported in Table \ref{Simulations1}.



\begin{table}[h!]\caption{Performance of fixed and adaptive designs }\label{Simulations1}
	\begin{center}
		\begin{tabular}{c|c|c||c|c||c}
			\multicolumn{3}{c||}{\bf{Parameters} }&\multicolumn{2}{c||}{\bf{MSE}($\widehat{\boldsymbol{\theta}}_{n}$)}&\multicolumn{1}{c}{\bf{Relative Efficiency}}\\
			\hline
			$\theta_2^t$& $\theta_{2,0}$ & $\sigma$ & Fixed & Adaptive & {Fix:Adap} \\
			\hline
			25 & 50 & 0.1 & 20.59& 20.61  & 0.999 \\
			\hline
			25 & 50 & 0.25 & 199.55& 398.89  & 0.500 \\
			\hline
			25 & 50 & 0.5 & 29094.20& 50970.56  & 0.571 \\
			\hline
			50 & 25 & 0.1 & 178.30& 158.95  & 1.122 \\
			\hline
			50 & 25 & 0.25 & 44837.38& 16634.79  & 2.69 \\
			\hline
			50 & 25 & 0.50 & 327364 & 210367.9 & 1.56 \\
			\hline
			10 & 25 & 0.1 & 2.72 & 2.38  & 1.14\\
			\hline
			10 & 25 & 0.25 & 18.77 & 38.08  & 0.49 \\
			\hline
			10 & 25 & 0.5 & 1323.66 & 8450.77  & 0.16 \\
			\hline
			25 & 10 & 0.1 & 30.09 & 23.93  & 1.26 \\
			\hline
			25 & 10 & 0.25 & 19106.01 & 321972  & 5.93 \\
			\hline
			25 & 10 & 0.5 & 254911.8 & 92494.97  & 2.76 \\
			\hline
		\end{tabular}
	\end{center}
	\label{confronti}
\end{table}

From the simulations, the adaptive design seems to perform better than the fixed one whenever the nominal value $\theta_{2,0}$ is inferior to the true value $\theta_2^t$. These results may be clarified by the following considerations. 

\begin{figure}
\includegraphics[width=10cm]{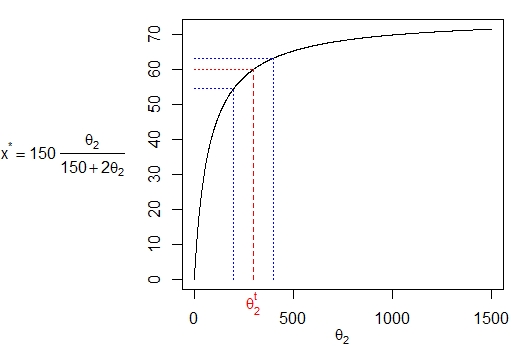}
\caption{D-optimum middle dose $x^*(\theta_2)$ for the Emax model}
\label{fig:Dopt}
\end{figure}

Note that, the derivative of $x^*(\theta_2)$ is a positive decreasing function of $\theta_2$ and thus the effect on $x^*(\theta_2)$ is larger  for the values $\theta_{2,0}<\theta_2^t$ (see Figure \ref{fig:Dopt}).
Moreover, the bias of the first stage MLE $\hat\theta_{n_1,2}$ is always positive as proved in Proposition \ref{prop:bias}. Hence, when $\theta_{2,0}<\theta_2^t$ the fixed procedure seems to have a worst performance, while $\hat\theta_{n_1,2}$ has a positive bias and thus it takes larger values and we expect that $x^*(\hat\theta_{n_1,2})$ is closer to $x^*(\theta_2^t)$ than $x^*(\theta_{2,0})$.

\begin{proposition}\label{prop:bias}
	If  $n_1$ first stage observations are taken according to  the D-optimal design (\ref{D-opt}) {with equal numbers treated at $a$, $x^*(\theta_{2,0})$ and  $b$}, then the bias of the first stage MLE of $\theta_2$ is
	$$
	E(\hat\theta_{n_1,2}-\theta_2)=\frac{b_{2}(\boldsymbol{\theta})}{n_1}+O(n_1^{-2}),
	$$ 
	where $b_{2}(\boldsymbol{\theta})>0$ given by
	\begin{eqnarray}
	b_{2}(\boldsymbol{\theta}) &=&
	\frac{1}{(a-b)^4 \theta_1^2 \theta_2^2 (a+\theta_{2,0})^2 (b+\theta_{2,0})^2} \nonumber\\
	&&\cdot 
	\big\{
	3\sigma^2 (a+\theta_2)^2 (b+\theta_2)^2 [2ab+(a+b)\theta_{2,0} +\theta_2(a+b+2\theta_{2,0})]^2 
	\nonumber\\
	&&\;
	[\,3ab (a+b)+(a^2+10 a b +b^2)\theta_{2,0}+3 (a+b) \theta_{2,0}^2\nonumber\\
	&& + 2 \theta_2 (a^2+ab+b^2+3(a+b)\theta_{2,0} +3 \theta_{2,0}^2)\,]\,
	\big\}.
	\label{bias}
	\end{eqnarray}
\end{proposition}
{\it Proof}.
Cox and Snell (1968) introduced the $O(n^{-1})$  formula  {for the bias of the} MLE in the case of $n$  observations {not  being identically distributed}. Cordeiro and Klein (1994) proposed a matrix expression for this bias, which is herein specialized for the Emax model and the D-optimal design $\xi_D^*(\theta_{2,0})$. Calculations are available by the authors upon request.


\section{Conclusions}\label{Sect:concl}
In this paper some important theoretical results about the maximum likelihood estimator are proved when observations are taken from a non linear gaussian model  according a two-stage procedure; the model involves a multidimensional parameter. The novelty from the previous literature is that the sample size at the first stage is small and thus standard asymptotic results cannot be applied.

First, the consistency of the MLE is proved under suitable assumptions, commonly satisfied. Then a central limit theorem is obtained, providing a closed form  of the asymptotic distribution of the MLE, which is a multivariate Gaussian mixture. As a corollary, the asymptotic covariance is found \emph{not} to be the inverse of the information matrix as in the standard cases, although they are connected in a specific way.

Finally, as an example, some simulations for the Emax dose-response model are performed to show how the method works. In this case, the exact expression of the first-order bias of the MLE at the first stage is also given. This result suggests, as a future development, some possible bias-corrections of the first stage estimate, that may hopefully improve the proposed two-stage procedure.


\bigskip

{\bf Acknowledgments} We are grateful to HaiYing Wang, Adam Lane and Giacomo Aletti for their precious comments which helped us to improve this work.

\bibliographystyle{abbrv}
\bibliography{ChiaraMayFlournoy}

\end{document}